\documentclass[12pt]{article}
\usepackage{amsmath,amsthm,amssymb}
\usepackage{color}
\newtheorem{theorem}{Theorem}[section]
\newtheorem{lemma}{Lemma}[section]
\newtheorem{proposition}{Proposition}[section]
\newtheorem{corollary}{Corollary}[section]
\newtheorem{remark}{Remark}[section]
\theoremstyle{definition}
\newtheorem{definition}{Definition}[section]
\newcommand{\beq}{\begin{eqnarray}}
\newcommand{\eeq}{\end{eqnarray}}
\newcommand{\beqq}{\begin{eqnarray*}}
\newcommand{\eeqq}{\end{eqnarray*}}
\usepackage[a4paper,top=3cm,bottom=3cm,left=2.5cm,right=2.5cm]{geometry}

\begin{document}
\mathsurround=1pt
\title{Hankel and Berezin type operators on weighted  Besov
spaces of
holomorphic functions on polydiscs}
\author{Anahit V. Harutyunyan\footnote{%
\textsc{Supported by DFG MA 2469/3-1}},\,\, George Marinescu}
\date{Yerevan State University,  University of Cologne}

\maketitle
\begin{abstract}{Assuming that  $S$ is the space of functions of
regular variation and  $\omega =(\omega_1,\ldots,\omega_n)$,
$\omega_j\in S$, by  $B_p(\omega )$ we denote the  class of all
holomorphic functions defined on the polydisk $U^n$ such that
$$\|f\|^p_{B_p(\omega
)}=\int_{U^n}|Df(z)|^p\prod_{j=1}^n\frac{\omega_j(1-|z_j|)dm_{2n}(z)}{(1-|z_j|^2)^{2-p}}<+\infty
,$$
where $dm_{2n}(z) $ is the $2n$-dimensional Lebesgue measure on $U^n$
and $D $ stands for a special fractional derivative  of $f$ defined here. 

\medskip\noindent In this paper we  consider  the generalized  little
Hankel  and Berezin type operators on $B_p(\omega)$ (and on
$L_p(\omega)$) and prove some theorems about the boundedness of these
operators. }
\vspace*{0.5cm}


\vspace*{0.5cm}
\end{abstract}
\section { Introduction  and auxiliary constructions}\label{s.1}
Numerous authors have contributed to holomorphic Besov spaces in the
unit disc in $\mathbb{C}$ and in the
unit ball  in $\mathbb{C^n}$, Arazy-Fisher-Peetre \cite{a2}, K. Stroethoff
\cite{st}
O. Blasco \cite{bl}, A. Karapetyants \cite{ka} see  K. Zhu \cite{zu}. The
investigation of holomorphic Besov space on the
polydisc is of special interest. The polydisc is a product of $n$ disks
and one would expect that
the natural generalisations of results from the one-dimensional case
would be valid here, but it
turns out that this is not true. The case of polydisc is different from
the $n=1$ case and from
the case of the $n$-dimensional ball. For example, let us
consider the classical theorem of Privalov: if $f\in$ {\rm Lip
}$\alpha$, then $Kf\in$ {\rm Lip }$\alpha,$
where $Kf$ is a Cauchy type integral. It is known that the analogue of
this theorem for
multidimensional Lipschitz classes is not true (\cite{yo}), even
though the analogue of this theorem for a sphere is valid (\cite{ru}).
In many cases, especially, when the class is defined by means of
derivatives,
  the generalisation of functional spaces in the polydisc is different
from that on a unit ball.
The generalisation of holomorphic Besov spaces on the polydisc see in
\cite{hl0}.
Let
$$
U^n=\{z=(z_1,\ldots,z_n)\in\mathbb{C}^n,\,\,|z_j|<1,\,\, 1\leq j\leq n\}
$$
be the unit polydisc in the $n$-dimensional complex plane $\mathbb{C^n}$ and
$$
T^n=\{z=(z_1,\ldots,z_n)\in\mathbb{C}^n,\,\, |z_i|=1,\,\,1\leq i\leq n\}
$$
be its torus. We denote by $H(U^n)$ the set of holomorphic functions on
$U^n$,
  by $L^\infty (U^n)$ the set of bounded  measurable functions on $U^n$
and by $H^\infty (U^n)$ the subspace of $L^\infty (U^n)$ consisting of
holomorphic functions.

Let $S$ be the class of all non-negative measurable  functions
$\omega $ on $(0,1),$ for which there exist positive numbers
$M_{\omega},\,\,q_{\omega},\,\,m_{\omega},\,\,(m_{\omega},q_{\omega}\in
(0,1)),$
such that
$$m_{\omega}\leq \frac {\omega (\lambda r)}{\omega (r)}\leq
M_{\omega},$$
for all $r\in (0,1)$ and $\lambda \in [q_{\omega}, 1].$ Some
properties of
functions from $S$ can be found in \cite{se}.
  We set
$$
-\alpha_{\omega}=\frac{\log m_{\omega}}{\log q_{\omega}^{-1}}; \qquad
\beta_{\omega}=
\frac{\log M_{\omega}}{\log q_{\omega}^{-1}}
$$
and assume that  $0<\beta_\omega <1.$
For example, $ \omega \in S$ if $\omega(t) = t^{ \alpha}$, where $ -1 <
\alpha <\infty$.

Using the results of \cite{se} one can prove that
$$
\omega_j(t)=\exp\biggl\{\eta_j(t)+\int_t^1\frac{\varepsilon_j(u)}{u}du\biggr\},
$$
where $\eta (u),$  $\varepsilon (u)$ are bounded measurable functions
and
$-\alpha_{\omega_j}\leq\varepsilon_j(u)\leq\beta_{\omega_j}\,\,(1\leq
j\leq n).$
Without loss of generality we assume that $\eta (u)=0.$  Then
\[ t^{ \alpha_{ \omega_j}} \leq \omega_j(t) \leq t^{- \beta_{
\omega_j}}\] is always true.

Below, for convenience of notations, for
$\zeta=(\zeta_1,...,\zeta_n),\,\,z=(z_1,...,z_n)$ we set
\[
\omega(1-|z|)=\prod_{j=1}^n\omega_j(1-|z_j|),\,\,
1-|z|=\prod_{j=1}^n(1-|z_j|),\,\,1-\overline\zeta
z=\prod_{j=1}^n(1-\overline \zeta_jz_j).
\]
Further, for
$m=(m_1,...,m_n)$ we set
\[
\begin{split}
&(m+1)=(m_1+1)...(m_n+1),\,\,
(m+1)!=(m_1+1)!...(m_n+1)!,\\
&(1-|z|)^m =
\prod_{j=1}^n(1-|z_j|)^{m_j}.
\end{split}
\]
Throughout the paper let assume $\omega_j\in S,\,1\leq j\leq n.$
The following definition gives the notion of the fractional
differential.

\begin{definition} {\it For a holomorphic function
$f(z)=\sum_{(k)=(0)}^{(\infty)}a_kz^k$, $z\in U^n$, and for $\beta =(\beta_1,...,\beta_n)$,
$\beta_j>-1$, $(1\leq j\leq n)$, we define the fractional  differential $D^\beta$  as  follows
$$D^\beta f(z)=\sum_{(k)=(0)}^{(\infty)}\prod_{j=1}^n
    \frac{\Gamma (\beta_j+1+k_j)}{\Gamma (\beta_j+1)\Gamma (k_j+1)}
a_kz^k,
    \,\,\, k=(k_1,...,k_n),\quad z\in U^n ,
    $$
    where $\Gamma (\cdot)$ is the  Gamma function and
$\sum_{(k)=(0)}^{(\infty)}=\sum_{k_1=0}^{\infty}\ldots\sum_{k_n=0}^{\infty}$}.
\end{definition}

 If  $\beta=(1,\ldots ,1)\,\,\mathrm{ then  \ we \ put }\,\, D^\beta
f(z)\equiv Df(z).$      Hence
      $$ Df(z_1, \ldots, z_n) = \frac{ \partial^n(f(z_1, \ldots, z_n) z_1
\cdots  z_n)}{ \partial z_1 \ldots \partial z_n}\,\cdot$$
If $n=1$ then $Df$ is the usual derivative of the function $zf(z)$.

    Let us define the weighted $L_p(\omega)$ spaces of holomorphic
functions.
    \begin{definition}{\it Let  $0<p<+\infty
1,\,\,\beta_{\omega_j}<-1(1\leq j\leq n).$ We denote by $L_p(\omega ) $
the set of all  measurable functions on $U^n$,  for which
    $$
    \|f\|^p_{L_p(\omega )}=\int_{U^n}
|f(z)|^p\frac{\omega(1-|z|)}{(1-|z|^2)^2}dm_{2n}(z)
      <+\infty .
    $$}

    \end{definition}

    Note that $L_p( \omega)$ is the $L_p-$space with respect to the
measure  \\
    $ \omega(1-|z|)(1-|z|^2)^{-2} dm_{2n}(z)$. Using the conditions
on $ \omega$ ($ \omega_j \in S $) we conclude that this measure is bounded.

    Now we define holomorphic Besov spaces on the polydisc.
    \begin{definition}{\it Let $0< p<+\infty$ and $f\in H(U^n).$ The
function $f$ is said to be in $B_p(\omega)$ if
    $$
    \|f\|^p_{B_p(\omega )}=\int_{U^n}
|Df(z)|^p\frac{\omega(1-|z|)}{(1-|z|^2)^{2-p}}dm_{2n}(z)
      <+\infty\,.
    $$}
    \end{definition}

    From the  definition of $Df$ it follows that $|| \cdot ||_{B_p(
\omega)}$ is indeed a norm. (We do not have to add $|f(0)|$). This
follows from the fact that here $Df =0$ implies  $f=0$
    for a holomorphic $f$.

    As in the one-dimensional case, $B_p(\omega)$ is a Banach space with
respect to the norm $\|\cdot\|_{B_p(\omega)}.$ For properties of
holomorphic Besov spaces see \cite{hl0}.

The investigation of Toeplitz operators  are widely known (see for
example \cite{h1,h2, st1, ma}).
Some problems on the Toeplitz operators can be solved by means of
Hankel operators and vice versa. In the classical theory of Hardy of
holomorphic functions on the unit disk there is only one type of Hankel
operator. In the $B_p(\omega)$ theory they are two:  little Hankel
operators and big  Hankel operators.  The analogue of the Hankel
operators of the Hardy theory here are little Hankel operators, which
were investigated by many authors (see for example  \cite{ja, af, hl0}).

Let us define the little Hankel operators as follows: denote by
$\overline B_p(\omega)$ the space of conjugate holomorphic functions on
$B_p(\omega).$ For the integrable function $f$ on $U^n$  we define the
generalized little Hankel operator  with symbol $h\in L^\infty(U^n) $ by
\beqq&&
h^\alpha_g(f)(z)=\overline P_\alpha (fg)(z)=\int_{U^n}\frac{(1-|\zeta
|^2)^\alpha}{(1-\zeta\overline z)^{\alpha +2}}f(\zeta)g(\zeta
)dm_{2n}(\zeta ),\\
&&
\alpha =(\alpha_1,\ldots ,\alpha_n),\,\alpha_j>-1,\,1\leq j\leq n.
\eeqq
For  $n=1,\,\alpha =0$ this includes the definition of the classical
little Hankel operator, see  \cite{z1}.
In Section 2 we consider the boundedness of little Hankel operator on
$B_p(\omega)$. For the case  $0<p<1$ and for the case $p=1$ we have the
following results
\begin{theorem}\label{th}Let $0<p< 1,\,\, f\in B_p(\omega)$ (or $f
\in\overline B^p(\omega)) ,\,g\in L^\infty (U^n).$  Then
$h^\alpha_g(f)\in \overline B_p(\omega)$  if and only if
$\alpha_j>\alpha_{\omega_j}/p-2,\,1\leq j\leq n.$\end{theorem}
\begin{theorem}\label{th12}Let $ f\in B_1(\omega),\, \,g\in L^\infty
(U^n).$  Then $h^\alpha_g(f)\in \overline B_1(\omega)$  if and only if
$\alpha_j>\alpha_{\omega_j}-2,\,1\leq j\leq n.$\end{theorem}
The case $p>1$ is different from the cases of $0>p<1$ and from the case
of $p=1$. Here we have the following
\begin{theorem} Let $1<p< +\infty,\,\, f\in B_p(\omega)$ (or $f
\in\overline B_p(\omega)) ,\,g\in L^\infty (U^n).$  Then  if
$\alpha_j>\alpha_{\omega_j},\,1\leq j\leq n$ then $h^\alpha_g(f)\in
\overline B_p(\omega)$. \end{theorem}

The Berezin transform is the analogue of the Poisson transform in the
$A^p(\alpha)$ (respectively, $( B_p(\omega))$) theory. It plays an
important role especially in the study of Hankel and Toeplitz operators.
In particular, some properties of those operators (for example,
compactness, boundedness) can be proved by means of the Berezin
transform (see \cite{ st, ma1,z1}). The Berezin-type operators,  on the
other hand, are of independent interest.

In the last Section 3 it will be shown, that some properties of
Berezin-type operators of the one dimensional classical case also hold
in the more general situation.  For the integrable function $f$  on
$U^n$  and for $g\in L^\infty (U^n)$ we define the Berezin-type operator
in the following way
$$B_g^\alpha f(z)=\frac{ (\alpha +1)}{\pi^n}(1-|z|^2)^{\alpha
  +2}\int_{U^n}\frac{(1-|\zeta |^2)^\alpha}{|1-z\overline
  \zeta|^{4+2\alpha}}f(\zeta)g(\zeta )dm_{2n}(\zeta).$$
In the case $\alpha =0,\,\,g\equiv 1$ the operator $B_g^\alpha $ will be called the
Berezin
transform.
We have the following results:

1. for the case of $0<p<1$ we have
\begin{theorem}Let $0<p< 1,\,\,f\in B_p(\omega)$ (or $f \in\overline
B_p(\omega )) ,\,g\in L^\infty (U^n)$ and let
$\alpha_j>\alpha_{\omega_j}/p-2,\,1\leq j\leq n.$  Then
$B^\alpha_g(f)\in L^p(\omega).$
\end{theorem}
2. the case  $1<p<+\infty$ gives the next theorem
  \begin{theorem}Let $1<p<+\infty ,\,\,f\in B_p(\omega)$ (or $f
\in\overline B_p(\omega )) ,\,g\in L^\infty (U^n)$ and let
$\alpha_j>(\alpha_{\omega_j}/p-2,\,1\leq j\leq n.$  Then
$B^\alpha_g(f)\in L_p(\omega).$
\end{theorem}

3. we consider now the case of $p=1.$
\begin{theorem}Let $f\in B_1(\omega)$ (or $f \in\overline B_1(\omega ))
,\,g\in L^\infty (U^n)$.  Then $B^\alpha_g(f)\in L_1(\omega)$ if and only
if    $\alpha_j>\alpha_{\omega_j},\,1\leq j\leq n.$
\end{theorem}

In general, $h^\alpha_g(f)$ and $B_g^\alpha $ are  not bounded.

To prove the main results we need some other notation.
The partition of the polydisc into dyadic quadrangles plays an
important role (see \cite{d2, sh}). Put
\begin{eqnarray*}
&&\Delta_{k_j,l_j}=\big\{z_j\in U: 1-\frac{1}{2^{k_j}}\leq |z_j|<
1-\frac{1}{2^{k_j+1}},\,\,\,
\frac{\pi l_j}{2^{k_j}}\leq\arg z_j<\frac{\pi (l_j +1)}{2^{k_j}}\bigr\},\\
&&
\Delta_{k_j,l_j}^*=4/3\Delta_{k_j,l_j},
\end{eqnarray*}
where
$k=(k_1,\ldots,k_n)$ $(k_j\geq 0)$, $l_j$ are some integers such that
$-2^{k_j}\leq l_j\leq2^{k_j+1}-1$ $(1\leq j\leq n)$ and
$2^k=(2^{k_1},\ldots,2^{k_n}).$

Then  $\quad\Delta_{k,l}=\Delta_{k_1,l_1}\times
\ldots\times\Delta_{k_n,l_n}
$
and $\Delta_{k,l}^*$ is defined similarly. The system $\{\Delta_{k,l}\}$ is called the system of dyadic quadrangles.

\begin{proposition}\label{p1} Let $\zeta_{k_j,l_j}$ be the center
  of  $\Delta_{k_j,l_j},\,\,1\leq j\leq n.$ Then
  $$1-|\zeta_{k_j,l_j}|\asymp 1-|\zeta_j |\quad \zeta_j\in
\Delta_{k_j,l_j}\quad \mathrm{and}\quad
  (1-|\zeta_{k_j,l_j}|)^2\asymp |\Delta_{k_j,l_j}|\quad 1\leq j\leq
n.$$\end{proposition}

Note that the partition of the polydisc into dyadic quadrangles is
important for obtaining some integral estimates particularly in the case
$0<p\leq 1$ \cite{sh}. Besides, the system $\{\Delta_{kl}\}$, as well as
the system $\{\Delta_{kl}^*\}$, are coverings of $U^n$, and one can
observe that the interiors of $\Delta_{kl}$ for disjoint indices are
disjoint, which is no longer true for  $\Delta_{kl}^*$. On the other
hand, $\{\Delta_{k,l}^*\}$ is a finite covering in the sense that any
quadrangle $\{\Delta_{kl}^*\}$ has nonempty intersection only with a
finite number of quadrangles from $\{\Delta_{kl}\}$, and this number is
independent of $k$ and $l$. Note that this partition for the spaces
$A^p_\alpha$ was used for the first time by F. A. Shamoyan \cite{sh} who
greatly investigated the theory of weighted classes of functions in the
polydisc and unit ball in $\mathbb{C}^n$.

    To prove the main results we need the following auxiliary lemmas:
    \begin{lemma} Let  $m=(m_1,\ldots,m_n)$ and
$\beta=(\beta_1,\ldots,\beta_n)$, $\beta_j\geq 0,1\leq j\leq n$. Then
If $f\in B_p(\omega)$ then
    \begin {equation}\label{e2}
    |f(z)|\leq C\int_{U^n}\frac{(1-|\zeta |^2)^m}{|1-\overline\zeta
z|^{m+1}}|D f(\zeta)|dm_{2n}(\zeta)
    \end {equation}
    where $m_j\geq \alpha_{\omega_j}-1\,\,(1\leq j\leq n).$
    \end{lemma}

The proof follows from  \cite[Lemma 2.5]{hl0}.

    \begin{lemma}
    Let $n=1$.  Assume $a+1-\beta_\omega >0,$ $b>1$ and
    $b-a-2>\alpha_\omega$. Then
        \begin {equation}\label{e3}\int_{U}\frac{(1-| \zeta|^2)^a\omega(1-|
\zeta|^2)}{|1-z\overline \zeta|^{b}}dm_2( \zeta)
    \leq\frac{\omega (1-|z|^2)}{(1-|z|^2)^{b-a-2}}\,\cdot
    \end {equation}
    \end{lemma}

    For the proof    see  \cite[Lemma 2]{h2}.

\section {Little Hankel  operators on $B_p(\omega)$ }\label{s.2}

  We consider the  little Hankel operators on  $B_p(\omega)(0<p<+\infty).$
  We denote the restriction of $ || \cdot ||_{L^p( \omega)}$ to
$ \overline{B}_p( \omega)$ by $|| \cdot ||_{ \overline{B}_p( \omega)}$.
First off all we consider the case $0<p<1$.
\begin{theorem}\label{th}Let $0<p< 1,\,\, f\in B_p(\omega)$ (or $f
\in\overline B^p(\omega)) ,\,g\in L^\infty (U^n).$  Then
$h^\alpha_g(f)\in \overline B_p(\omega)$  if and only if
$\alpha_j>\alpha_{\omega_j}/p-2,\,1\leq j\leq n.$\end{theorem}
{\bf Proof.} Let $0<p< 1,\,\,f\in B^p(\omega)$ (or $f \in\overline
B^p(\omega)),\,\,g\in L^\infty(\omega )$ and
$\alpha_j>\alpha_{\omega_j}/p-2,\,1\leq j\leq n.$
We will show that $h^\alpha_g (f)\in \overline B^p(\omega)$. Using the
partition of the polydisc,  Lemma 3  from \cite{sh}  and Proposition
\ref{p1}, we get
\beqq&&
I=\int_{U^n}\frac{\omega
(1-|z|)}{(1-|z|^2)^{2-p}}\biggl(\int_{U^n}\frac{(1-|\zeta
|^2)^\alpha}{|1-\overline z
  \zeta|^{\alpha +3}}|f(\zeta)||g(\zeta
)|dm_{2n}(\zeta)\biggr)^pdm_{2n}(z)\leq\\
&& C(g) \int_{U^n}\frac{\omega
(1-|z|)}{(1-|z|^2)^{2-p}}\sum_{k,l}\biggl(\int_{\Delta_{k,l}}\frac{(1-|\zeta
  |)^\alpha}{|1-\overline z\zeta |^{\alpha +2}}|f(\zeta )|dm_{2n}(\zeta
)\biggr)^pdm_{2n}(z)\leq\\
&& C(g)
\int_{U^n}\frac{\omega
(1-|z|)}{(1-|z|^2)^{2-p}}\sum_{k,l}\max_{\substack{\zeta\in
\overline\Delta_{k,l} }}|f(\zeta
)|^p|\Delta_{k,l}|^p\frac{(1-|\zeta_{k,l}
  |)^{\alpha p}}{|1-\overline z\zeta_{k,l} |^{(\alpha
+3)p}}dm_{2n}(z)=\\
 && C(g)\sum_{k,l}\max_{\substack{\zeta\in \overline\Delta_{k,l}
}}|f(\zeta
    )|^p|\Delta_{k,l}|^p(1-|\zeta_{k,l}
        |)^{\alpha p}\int_{U^n}\frac{\omega
(1-|z|)(1-|z|^2)^{p-2}}{|1-\overline z\zeta_{k,l} |^{(\alpha
+3)p}}dm_{2n}(z)
  \eeqq
where  $\zeta_{k,l}$ is the center of $\Delta _{k,l}$ and
$I=\|h^\alpha_{g}f\|_{\overline B^p(\omega)},\,\, C(\alpha, p,
\omega)\|g\|_\infty=C(g)$

Recalling that $\{\Delta^*_{k,l}\}$ forms a finite covering of $U^n$, by
(\ref{e3}) and Lemma  4 from \cite{sh} we obtain
\beqq&&
I\leq C(g)\sum_{k,l}\max_{\substack{\zeta\in\overline\Delta_{k,l}
}}|f(\zeta
)|^p|(1-|\zeta_{k,l}|)^{-2+2}\omega(1-|\zeta_{k,l})\leq\\
&& C(g)\sum_k\sum_l\int_{\Delta^*_{k,l}}|f(z)|^p \frac{\omega
(1-|z|)}{(1-|z|^2)^{2}})dm_{2n}(\zeta)\leq\\
&& C(g)\int_{U^n}|f(z)|^p \frac{\omega
(1-|z|)}{(1-|z|^2)^{2}})dm_{2n}(\zeta)
\eeqq
Using (\ref{e2}) we get
\beqq&&
I\leq C(g) \int_{U^n} \frac{\omega
(1-|z|)}{(1-|z|^2)^{2}}\Biggl(\int_{U^n}
\frac{(1-|t|^2)^m}{|1-\overline t\zeta|^{m+1}}|Df(t)|dm_{2n}(t)
\Biggr)^pdm_{2n}(\zeta)\leq\\
&& C(g) \int_{U^n} \frac{\omega
(1-|z|)}{(1-|z|^2)^{2}}\sum_{k,l}\Biggl(\int_{\Delta_{k,l}}\frac{(1-|t|^2)^m}{|1-\overline
t\zeta|^{m+1}}|Df(t)|dm_{2n}(t)\Biggr)^pdm_{2n}(\zeta)\leq \\
&& C(g) \int_{U^n} \frac{\omega (1-|z|)}{(1-|z|^2)^{2}}\sum_{k,l}\max_{t\in
\overline{\Delta_{k,l}}}|Df(t)|^p|\Delta_{k,l}|^p\frac{(1-|t_{k,l}|^2)^{mp}}{|1-\overline
t\zeta|^{(m+1)p}}dm_{2n}(\zeta)\leq\\
&& C(g) \sum_{k,l}\max_{t\in
\overline{\Delta_{k,l}}}|Df(t)|^p|\Delta_{k,l}|^p \frac{\omega
(1-|t_{k,l}|)(1-|t_{k,l}|^2)^{mp}}{(1-|t_{k,l}2)^{(m+1)p-2+2}}=\\
&& C(g) \sum_{k,l}\max_{t\in
\overline{\Delta_{k,l}}}|Df(t)|^p(1-|t_{k,l}|^2)^{p-2+2}.
\eeqq
In the last inequality we have used Lemma 2  again.
By Lemma 4  from \cite{sh} we get
\beqq&&
I\leq  C(g))\sum_k\sum_l\int_{\Delta^*_{k,l}}\frac{\omega
(1-|z|)}{(1-|z|^2)^{2-p}}|Df(z)|dm_{2n}(z)\leq\\
 && \int_{U^n} \frac{\omega (1-|z|)}{(1-|z|^2)^{2-p}}|Df(z)|dm_{2n}(z)=
C(\alpha, p, \omega)\|g\|_\infty\|f\|_{B_p(\omega)}^p,
\eeqq
which proves our statement.

Conversely, let $h^\alpha_g(f)\in \overline B_p(\omega)$ for all $g\in
L^\infty (U^n)$.
For $r=(r_1,...,r_n),\,\,r_j\in (0,1),\,\,k=(k_1,...,k_n)$  we take
the  function
\begin{equation}\label{ef1}f_r(z)=C_r(1-rz)^{-k},\,\,
k_j>(\alpha_{\omega_j}+2)/p,\,\, 1\leq j\leq n,
\end{equation}
where $C_r=(1-r)^{k}\omega^{-1/p}(1-r).$
Then  we have $\|f_r\|_{ B_p(\omega)}\thicksim const.$

We consider the following domains
$$
\widetilde U_j=\{z_j\in U, |\arg
z_j|<(1-r_j)/2;\,(4r_j-1)/3<|z_j|<(1+2r_j)/3\}
$$
and
$$
\widetilde U^n=\widetilde U_1\times\ldots \times \widetilde U_n.
$$
  Take  the function $g_r(\zeta )$ as
$$g_r(\zeta )=\exp^{-\arg f_r(\zeta )}
$$
and a polydisc $V^n$ centered at $(r_1,...,r_n)$ with radius of $(1-r_1)...(1-r_n)$
such that $\overline V^n\subset \widetilde U^n$( $\overline V^n$ is the
closure of $V^n$),  we get
$$
\|h^\alpha_{g_r}f_r\|_{\overline B_p(\omega)}\geq
\int_{U^n} \frac{\omega (1-|z|)}{(1-|z|^2)^{2-p}}\biggl(\int_{
V^n}\frac{(1-|\zeta |)^\alpha}{|1-\overline z\zeta |^{\alpha
+3}}|f_r(\zeta
)|dm_{2n}(\zeta)\biggr)^{p}dm_{2n}(z).
$$
  Let
$$\max_{\substack\zeta\in\overline V^n}|1-\overline
z\zeta|=|1-\overline z\widetilde \zeta |,
$$ then
\[
\begin{split}
\|h^\alpha_{g_r}f_r\|_{\overline B_p(\omega)}&\geq C_1(\alpha, p,
\omega)\frac{(1-r)^{\alpha p}}{\omega (1-r)}\int_{U^n}\frac{\omega
(1-|z|)}{(1-|z|)^{2-p}}\biggl(\int_{
V^n}\frac{dm_{2n}(\zeta)}{|1-\overline z\zeta |^{\alpha
+3}}\biggr)^pdm_{2n}(z)\\
&\geq
C_1(\alpha, p, \omega)\frac{(1-r)^{(\alpha +2)p}}{\omega
(1-r)}\int_{U^n}\frac{\omega (1-|z|)dm_{2n}(z)}{|1-\overline
z\widetilde\zeta|^{(\alpha +3) p}(1-|z|)^{2-p}}\,\cdot
\end{split}
\]
If we assume that $(\alpha_j+2)p\leq\alpha_{\omega_j}$ for some $j$,
then  for the corresponding integral taking
$\omega_j(t)=t^{\alpha_{\omega_j}}$ we get

$$
\int_{U^n}\frac{\omega (1-|z|)dm_{2n}(z)}{|1-\overline
z\widetilde\zeta|^{(\alpha +3) p}(1-|z|)^{2-p}}\thicksim const,\quad
\mathrm{if}\quad
(\alpha_j+2)p<\alpha_{\omega_j}
$$
and
$$
\int_{U^n}\frac{\omega (1-|z|)dm_{2n}(z)}{|1-\overline
z\widetilde\zeta|^{(\alpha +3) p}(1-|z|)^{2-p}})\thicksim
\log\frac{1}{1-|\widetilde\zeta_j|},\quad \mathrm{if}\quad
(\alpha_j+2)p=\alpha_{\omega_j}+2.
$$
Consequently,
$$
\frac{(1-r_j)^{(\alpha_j+2) p}}{\omega_j (1-r_j)}\rightarrow
\infty,\quad \frac{(1-r_j)^{(\alpha_j+2) p}}{\omega_j
(1-r_j)}\log\frac{1}{1-r_j}\rightarrow \infty
$$
if $r_j\rightarrow 1-0 .$

\hfill$\square$

\begin{corollary} Let $ 0<p< 1,$
$\alpha_j>\alpha_{\omega_j}/p-2,\,1\leq j\leq n, \,\,g\in L^\infty
(U^n).$ Then $h^\alpha_g$ is bounded on $B_p(\omega)$, (and on
$\overline B_p(\omega)).$ Moreover, $\|h_g^\alpha\|\leq C
\|f\|\cdot\|g\|$
\end{corollary}

In the case if $p=1$ we have
\begin{theorem}\label{th12}Let $ f\in B_1(\omega),\, \,g\in L^\infty
(U^n).$  Then $h^\alpha_g(f)\in \overline B_1(\omega)$  if and only if
$\alpha_j>\alpha_{\omega_j}-2,\,1\leq j\leq n.$\end{theorem}

{\bf Proof.} Let $f\in B_1(\omega),\,g\in L^\infty (U^n)$ and
$C(\alpha, \omega)\|g\|_\infty=\widetilde C$. Then by (\ref{e2}) and
(\ref{e3}) we have
\beqq&&
\|h^\alpha_g(f)\|_{\overline B_1(\omega)}\leq
\|g\|_\infty\int_{U^n}(1-|\zeta |^2)^\alpha|f(\zeta
)|\int_{U^n}\frac{\omega(1-|z|)dm_{2n}(z)}{|1-\zeta\overline z|^{\alpha
+3}(1-|z|)}dm_{2n}(\zeta)\leq\\
&&\widetilde C\int_{U^n}|f(\zeta )|\frac{\omega (1-|\zeta |^2)}{(1-|\zeta
|)^2}dm_{2n}(\zeta)\leq\widetilde C\int_{U^n}\frac{\omega (1-|\zeta
|)}{(1-|\zeta |)^2}\int_{U^n}
\frac{(1-|t|^2)^m}{|1-\overline
t\zeta|^{m+1}}|Df(t)|\times\\
&&dm_{2n}(t)
dm_{2n}(\zeta)=\widetilde C\int_{U^n}
(1-|t|^2)^m|Df(t)|\int_{U^n}\frac{\omega (1-|\zeta
|)dm_{2n}(\zeta)dm_{2n}(t)}{|1-\overline t\zeta|^{m+1}(1-|\zeta
|^2)^2}.
\eeqq
Using (\ref{e2}) again we get
\begin{eqnarray}
\|h^\alpha_g(f)\|_{\overline B_1(\omega)}\leq\widetilde C\int_{U^n}
\frac{\omega(1-|t|)}{(1-|t|}|Df(t)|dm_{2n}(t)=\widetilde
C\|f\|_{\overline B_1(\omega)}.\nonumber
\end{eqnarray}
Next, assume that $h^\alpha_gf\in \overline B_1(\omega).$ The proof of
the necessity of the condition $\alpha_j>\alpha_{\omega_j},\,1\leq j\leq
n$ is similar to the corresponding proof in Theorem \ref{th}. We omit
the details. This proves the  theorem.

\hfill $\square$
\begin{corollary} Let  $\alpha_j>\alpha_{\omega_j},\,1\leq j\leq n,
\,\,g\in L^\infty (U^n).$ Then $h^\alpha_g$ is bounded on $B_1(\omega)$
and $\|h_g^\alpha\|\leq C \|f\|\cdot\|g\|$.
\end{corollary}

Now we consider the case of  $1<p< +\infty$.
\begin{theorem} Let $1<p< +\infty,\,\, f\in B_p(\omega)$ (or $f
\in\overline B_p(\omega)) ,\,g\in L^\infty (U^n).$  Then  if
$\alpha_j>\alpha_{\omega_j},\,1\leq j\leq n$ then $h^\alpha_g(f)\in
\overline B_p(\omega)$ \end{theorem}
{\bf Proof.} Let $1<p< +\infty,\,\, f\in B_p(\omega)$ (or $f
\in\overline B_p(\omega)) ,\,g\in L^\infty (U^n).$ We show that
$h^\alpha_g(f)\in \overline B_p(\omega)$.
By H\"older inequality and by (2)  we get
\beqq &&
|Dh_g^\alpha(f)(z)|\le\int_{U^n}\frac{(1-|\xi|^2)^\alpha}{|1-\xi\bar{z}|^{\alpha+3}}|f(\xi)|\cdot
|g(\xi)|\, dm_{2n}(\xi)\le\\
&&\Vert g\Vert_{\infty}\int_{U^n}\frac{(1-|xi|^2)^\alpha
|f(\xi)|}{|1-\xi\bar{z}|^{\alpha+3}}\, dm_{2n}(\xi)\le
\Vert g\Vert_\infty \times\\
&& \left(\int_{U^n}\frac{(1-|\xi|^2)^\alpha
|f(\xi)|^p}{|1-\xi\bar{z}|^{\alpha+3}}\, dm_{2n}(\xi)\right)^{1/p}\cdot
\left(\int_{U^n}\frac{(1-|\xi|^2)^\alpha\,
dm_{2n}(\xi)}{|1-\xi\bar{z}|^{\alpha+3}}\right)^{1/q}\le\\
&&\frac{C(\alpha,q)\Vert g\Vert_\infty}{(1-|z|)^{1/q}}\,
\left(\int_{U^n}\frac{(1-|\xi|^2)^\alpha
|f(\xi)|^p}{|1-\xi\bar{z}|^{\alpha+3}}\, dm_{2n}(\xi)\right)^{1/p}
\eeqq
Then setting $C(\alpha,q)\Vert g\Vert_\infty=C$ we have
\beqq&& \Vert
h_g^{(\alpha)}(f)\Vert_{B_p(\omega)}=\int_{U_n}\frac{\omega(1-|z|)
}{(1-|z|^2)^{2-p}}|D h_g^\alpha(f)(z)|^p\, dm_{2n}(z)\le\\
&&C \int_{U^n}\frac{\omega(1-|z|)}{|1-|z|^2)^{2-p+p/q}}\,
\int_{U^n}\frac{(1-|\xi|^2)^\alpha
|f(\xi)|^p}{|1-\xi\bar{z}|^{\alpha+3}}\, dm_{2n}(\xi) dm_{2n}(z)\le\\
&&C\int_{U^n}|f(\xi)|^p(1-|\xi|^2)^\alpha\int_{U^n}\,\frac{\omega(1-|z|)
dm_{2n}(z)
dm_{2n}(\xi)}{|1-\xi\bar{z}|^{\alpha+3}(|1-|z|^2)^{2-p+p/q}}\le\\
&& C_1\int_{U^n} (1-|\xi|^2)^\alpha
|f(\xi)|^p\frac{\omega(1-|\xi|)(1-|\xi|^2)^{p-2-p/q}}{(1-|\xi|^2)^{\alpha+1}}\,
dm_{2n}(\xi)=\\
&&\int_{U^n} (1-|xi|^2)^{p-2-p/q-1}\omega(1-|\xi|) |f(\xi)|^p\,
dm_{2n}(\xi).
\eeqq
In the last inequality we have used (\ref{e2}). On the other hand, by
(\ref{e2}) we get
\beqq &&|f(\xi)|^p\le
\left(\int_{U^n}\frac{(1-|t|^2)^m}{|1-\bar{t}\xi|^{m+1}}\, |D f(t)|\,
dm_{2n}(t)\right)^p\le\\
&&\left(\int_{U^n}\frac{(1-|t|^2)^{m-\delta}
(1-|t|^2)^\delta}{|1-\bar{t}\xi|^{m+1}} |D f(t)|\,
dm_{2n}(t)\right)^p\le\\
&&\int_{U^n}\frac{(1-|t|^2)^{m-\delta}}{|1-\bar{t}\xi|^{m+1}}
(1-|t|^2)^{\delta p} |D f(t)|^p\, dm_{2n}(t)\cdot
\frac{C(m,\delta,q)}{(1-|\xi|^2)^{(\delta-1)p/q}},
\eeqq
for some $\delta >1$.
Then we obtain
\beqq && \Vert h_g^{\alpha)} (f)\Vert{B_p(\omega)}\le
C_1\int{U^n}(1-|t|^2)^{m-\delta+\delta p} |D f(t)|^p\int_{U^n}
\frac{(1-|\xi|^2)^{p-3-\delta p/q}}{|1-\bar{t}\xi|^{m+1}}
\times\\&&\omega(1-|\xi|)dm_{2n}(\xi)
\int_{U^n} (1-|t|^2)^{m-\delta+\delta p} |D
f(t)|^p\le\\&&C_2\int_{U^n}
\frac{\omega(1-|\xi|)(1-|\xi|)^{(1-\delta)p/q-2}}{|1-\bar{t}\xi|^{m+1}}
dm_{2n}(\xi) d_{2n}(t)\le\\
&&\int_{U^n} (1-|t|^2)^{m-\delta+\delta p}\frac{|Df(t)|^p\omega(1-|t|)
dm_{2n}(t)}{(1-|t|^2)^{m-1+2-(1-\delta)p/q}}=\\
&&\int_{U^n} (1-|t|^2)^{p-2}|Df(t)|^p\omega(1-|t|)\, dm_{2n}(t)=\Vert
f\Vert_{B_p}(\omega).
\eeqq

We have $\Vert h_g^\alpha(f)\Vert_{B_p(\omega)}\le C_3\Vert
f\Vert_{B_p(\omega)}\|g\|_\infty$, where $C_3=C_2\cdot C^p(m,\delta,
q)$.

\hfill
$\square$

\begin{corollary} Let  $\alpha_j>\alpha_{\omega_j},\,1\leq j\leq n,
\,\,g\in L^\infty (U^n).$ Then $h^\alpha_g$ is bounded on $B_p(\omega)$
and $\|h_g^\alpha\|t_{B_p(\omega)}\leq C_3
\|f\|_{B_p(\omega)}\cdot\|g\|_{\infty}$
\end{corollary}

\section {Berezin-type  operators on $B_p(\omega)$ }\label{s.3}

In this section we consider the boundedness of the Berezin-type
operators. Let us  consider first the case  $0<p<1.$
\begin{theorem}\label{t1}Let $0<p< 1,\,\,f\in B_p(\omega)$ (or $f
\in\overline B_p(\omega )) ,\,g\in L^\infty (U^n)$ and let
$\alpha_j>\alpha_{\omega_j}/p-2,\,1\leq j\leq n.$  Then
$B^\alpha_g(f)\in L^p(\omega).$
\end{theorem}
{\bf Proof.} Let $f\in B_p(\omega)$ or $f\in \overline B_p(\omega).$ We
will show that $B_\alpha f\in
L_p(\omega ).$ To this end we estimate the corresponding integral
$$
\int_{U^n}\frac{\omega (1-|z|)}{(1-|z|^2)^{2}}\biggl((1-|z|^2)^{\alpha
  +2}\int_{U^n}\frac{(1-|\zeta
|^2)^\alpha|f(\zeta)||g(\zeta)|}{|1-z\overline
  \zeta|^{4+2\alpha}}dm_{2n}(\zeta)\biggr)^pdm_{2n}(z)\equiv I
$$
Using the  partition of the polydisc, we obtain
\beqq&&
I\leq \|g\|_\infty\int_{U^n}(1-|z|^2)^{(\alpha
+2)p-2}\omega(1-|z|)\times\\
&&\sum_{k,l}\biggl(\int_{\overline\Delta_{k,l}}\frac{(1-|\zeta
  |)^\alpha}{|1-\overline\zeta z|^{4+2\alpha}}|f(\zeta )|dm_{2n}(\zeta
)\biggr)^pdm_{2n}(z)\leq C(\alpha, \omega,
p)\|g\|_{\infty}\times\\
&&\int_{U^n}(1-|z|^2)^{(\alpha +2)p-2}\omega
(1-|z|)
\sum_{k,l}\max_{\substack{\zeta\in\overline\Delta_{k,l} }}|f(\zeta
)|^p|\Delta_{k,l}|^p\frac{(1-|\zeta_{k,l}
  |)^{\alpha p}dm_{2n}(z)}{|1-\overline\zeta_{k,l}
z|^{(4+2\alpha)p}}=\\
 && C(\alpha, \omega, p)\|g\|_{\infty}
  \sum_{k,l}\max_{\substack{\zeta\in\overline\Delta_{k,l} }}|f(\zeta
  )|^p|\Delta_{k,l}|^p\times\\
  &&\int_{U^n}(1-|z|^2)^{(\alpha
+2)p-2}\omega
    (1-|z|)\frac{(1-|\zeta_{k,l}
    |)^{\alpha p}dm_{2n}(z)}{|1-\overline\zeta_{k,l}
z|^{(4+2\alpha)p}}\,\cdot
  \eeqq

\medskip\noindent Taking into account that
$p(4+2\alpha_j)>(\alpha_j+2)p+\alpha_{\omega_j}\,\, (1\leq j\leq n)$
and (\ref{e2}),
we get
\beqq&&
I\leq C(\alpha ,
p,\omega)\|g\|_{\infty}\sum_{k,l}\max_{\substack{\zeta\in\overline\Delta_{k,l}
}}|f(\zeta
)|^p|(1-|\zeta_{k,l}|)^{2-2}\omega(1-|\zeta_{k,l})\leq\\
&&C(\omega, \alpha , p)\|g\|_{\infty}\int_{U^n}|f(z)|^p\frac{\omega
(1-|\zeta |)}{(1-|z|^2)}dm_{2n}(\zeta).
\eeqq

In the last inequality we have used Lemma 4 \cite{sh}. Next we estimate
the last integral. Using Lemma 1  we obtain
$$
I\leq C(\omega, \alpha , p)\|g\|_{\infty}\int_{U^n}\frac{\omega
(1-|\zeta
|)}{(1-\zeta|^2)}\biggl(\int_{U^n}\frac{(1-|t|^2)^m}{|1-\overline
{t}\zeta|^{m+1}}|Df(t)|dm_{2n}(t)\biggr)^pdm_{2n}(\zeta).
$$
Then from
\beqq&&
\biggl(\int_{U^n}\frac{(1-|t|^2)^m}{|1-\overline
{t}\zeta|^{m+1}}|Df(t)|dm_{2n}(t)\biggr)^p\leq \\
&&\sum_{k,l}\biggl(\int_{\Delta_{k,l}}\frac{(1-|t|^2)^m}{|1-\overline
{t}\zeta|^{m+1}}|Df(t)|dm_{2n}(t)\biggr)^p\leq\\
&&\sum_{k,l}\max_{\substack{t\in\overline\Delta_{k,l}
}}|Df(t)|^p|\Delta_{k,l}|^p\frac{(1-|t_{k,l}|^2)^{mp}}{|1-\overline
{t_{k,l}}\zeta|^{(m+1)p}},
\eeqq
we conclude
\beqq
I\leq C(\omega,
\alpha,p)\|g\|_{\infty}\sum_{k,l}\max_{\substack{t\in\overline\Delta_{k,l}
}}|Df(t)|^p|\Delta_{k,l}|^p(1-|t_{k,l}|^2)^{mp}\times\\
\int_{U^n}\frac{\omega(1-|z|)(1-|\zeta|)^{-2}}{|1-\overline
{t_{k,l}}\zeta|^{(m+1)p}}dm_{2n}(z)\leq\\ C(\omega,
\alpha,p)\|g\|_{\infty}
\sum_{k,l}\max_{\substack{t\in\overline\Delta_{k,l}
}}|Df(t)|^p|\Delta_{k,l}|^p(1-|t_{k,l}|^2)^{mp}\frac{\omega(1-|t_{k,l}|)}{(1-|t_{k,l}|^2)^{(m+1)p}}=\\
C(\omega, \alpha,p)\|g\|_{\infty}
\sum_{k,l}\max_{\substack{t\in\overline\Delta_{k,l}
}}|Df(t)|^p\omega(1-|t_{k,l}|)(1-|t_{k,l}|^2)^{p}\leq\\
\int_{U^n}|Df(t)|^p\frac{\omega(1-|t|)}{(1-|t|^2)^{2-p}}dm_{2n}(t)
\eeqq
In the last inequality we have used Lemma 4 again. Next we have
$$I\leq  C(\omega, \alpha,p)\|g\|_{\infty} \|f\|_{B_p(\omega)}$$

\hfill
$\square$

\begin{remark} The condition $\alpha_j+2>(\alpha_{\omega_j}+2)/p,
(1\leq j\leq n)$ in Theorem \ref{t1} is necessary too. Moreover, if
$B_\alpha $ is bounded on $L^p(\omega)$ then
$\alpha_j+2>(\alpha_{\omega_j}+2)/p,
(1\leq j\leq n).$
\end{remark}
\medskip\noindent The proof is similar to the corresponding part of
Theorem \ref{th} and we omit it.

\begin{corollary} Let $ 0<p< 1,$
$\alpha_j>(\alpha_{\omega_j}+2)/p-2,\,1\leq j\leq n, \,\,g\in L^\infty
(U^n).$ Then $B^\alpha_g$ is bounded on $A^p(\omega)$ and on $\overline
A^p(\omega).$
\end{corollary}

\begin{theorem} Let $1<p<+\infty ,\,\,f\in B_p(\omega)$ (or $f
\in\overline B_p(\omega )) ,\,g\in L^\infty (U^n)$ and let
$\alpha_j>(\alpha_{\omega_j}/p-2,\,1\leq j\leq n.$  Then
$B^\alpha_g(f)\in L_p(\omega).$
\end{theorem}
{\bf Proof.} Let $f\in B_p(\omega)$ or $f\in \overline  B_p(\omega).$
Our aim is to show that $B_\alpha f\in
L^p(\omega ).$ We have
\beqq && |B_g^\alpha(f)(z)|^p\le(1-|z|^2)^{(\alpha+2)p}\frac{C(\alpha,
\pi, p)}{(1-|z|^2)^{(\alpha+2)p/q}}\le\\
&&\int_{U^n}\frac{(1-|\xi|^2)^\alpha|f(\xi)|^p
|g(\xi)|^p}{|1-z\bar{\xi}|^{2\alpha+4}}\, dm_{2n}(\xi)\le C(\alpha, \pi,
p)(1-|z|^2)^{\alpha+2}\times\\
&&\int_{U^n}\frac{(1-|\xi|^2)^\alpha
|f(\xi)|^p}{|1-\bar{\xi}z|^{2\alpha+4}}\, dm_{2n}(\xi)\le C(\alpha, \pi,
p)(1-|z|^2)^{\alpha+2}\cdot \Vert g\Vert_\infty\times\\
&&\int_{U^n}\frac{(1-|\xi|^2)^\alpha}{|1-z\bar{\xi}|^{2\alpha+4}}\,
\int_{U^n}\frac{(1-|t|^2)^{m-\delta}(1-|t|^2)^{\delta p} |D
f(t)|^p}{|1-\bar{t}\xi|^{m+1}(1-|\xi|^2)^{(\delta-1)p/q}}\,
dm_{2n}(t)dm_{2n}(\xi)=\\
&&C(\alpha, \pi, p) (1-|z|^2)^{\alpha+2} \Vert g\Vert_\infty\int_{U^n}
(1-|t|^2)^{m-\delta+\delta p} |Df(t)|^p\times\\
&&\int_{U^n}\frac{(1-|\xi|^2)^{\alpha-(\delta-1)p/q}}{|1-\bar{\xi}z|^{2\alpha+4}|1-\bar{t}\xi|^{m+1}}\,
dm_{2n}(\xi) dm_{2n}(t).
\eeqq
Then
\beqq&&\Vert
B_g^\alpha(f)\Vert_{L_p(\omega)}=\int_{U^n}(1-|t|^2)^{m-\delta+\delta p}
|Df(t)|^p\int_{U^n}\frac{(1-|\xi|^2)^{\alpha-(\delta-1)p/q}}{|1-\bar{t}\xi|^{m+1}}\times\\
&&\int_{U^n}\frac{\omega(1-|z|)(1-|z|^2)^\alpha}{|1-\bar{\xi}z|^{2\alpha+4}}\,
dm_{2n}(z)dm_{2n}(\xi)dm_{2n}(t)\le
\int_{U^n}(1-|t|^2)^{m-\delta+\delta p}\times\\&&|Df(t)|^p
\int_{U^n}\frac{(1-|xi|^2)^{\alpha-(\delta-1)p/q}\omega(1-|\xi|)}{|1-\bar{t}\xi|^{m+1}(1-|\xi|)^{\alpha+2}}\,
dm_{2n}(\xi)dm_{2n}(t)=\\
&&\int{U^n}(1-|t|^2)^{m-\delta+\delta p}
|Df(t)|^p\int_{U^n}\frac{\omega(1-|\xi|)}{|1-\bar{t}\xi|^{m+1}}(1-|xi|^2)^{-2-(\delta-1)p/q}\,
dm_{2n}(\xi)\\
&&\int_{U^n}(1-|t|^2)^{m-\delta+\delta p}\frac{\omega(1-|t|)
|Df(t)|^p}{(1-|t|^2)^{m-1+2(\delta-1)p/q}}\, dm_{2n}(t)=\\
&&\int_{U^n}\frac{\omega(1-|t|) |Df(t)|^p}{(1-|t|^2)^{2-p}}\,
dm_{2n}(t)=\Vert f\Vert_{B_p(\omega)}\Vert g\Vert_\infty C(\alpha, \pi,
p).
\eeqq

\hfill$\square$

We consider now the case of $p=1.$
\begin{theorem}Let $f\in B_1(\omega)$ (or $f \in\overline B_1(\omega ))
,\,g\in L^\infty (U^n)$.  Then $B^\alpha_g(f)\in L_1(\omega)$ if and only
if    $\alpha_j>\alpha_{\omega_j},\,1\leq j\leq n.$
\end{theorem}
{\bf Proof.} Let $f\in B_1(\omega)$ or $f\in \overline B_1(\omega).$
Our aim is to show that $B^\alpha_g f\in
L_1(\omega ) .$ We have
\beqq &&
|B^\alpha_g(f)|\le(1-|z|^2)^{\alpha+2}\int_{U^n}\int_{U^n}\frac{(1-|xi|^2)^\alpha
|f(\xi)|\cdot |g(\xi)|\, dm_{2n}(\xi)}{|1-\bar{xi}|z|^{4+2\alpha}}\le\\
&&\Vert g\Vert_\infty
(1-|z|^2)^{\alpha+2}\int_{U^n}\frac{(1-|\xi|^2)^\alpha}{|1-\bar{\xi}z|^{4+2\alpha}}
\int_{U^n}\frac{(1-|t|^2)^m |Df(t)|}{|1-\bar{t}\xi|^{m+1}}\, dm_{2n}(t)
dm_{2n}(\xi).
\eeqq
Then using  \eqref{e3} we get
\[
\begin{split}
\Vert
B_g^\alpha\Vert_{L_1(\omega)}&\le\int_{U^n}
\int_{U^n}\frac{(1-|\xi|^2)^\alpha}{|1-\bar{\xi}z|^{2\alpha+4}}\int_{U^n}\frac{(1-|t|^2)^m
|Df(t)| \omega(1-|z|)}{|1-\bar{t}\xi|^{m+1} (1-|z|^2)^{2-\alpha-2}}\,
dm_{2n}(z) dm_{2n}(t) dm_{2n}(\xi) \\
&=\int_{U^n}(1-|t|^2)^m
|Df(t)|\int_{U^n}\frac{(1-|\xi|^2)^\alpha}{|1-\bar{t}\xi|^{m+1}}
\int_{U^n}\frac{\omega(1-|z|)\, dm_{2n}(z)
dm_{2n}(\xi)dm_{2n}(t)}{|1-\bar{\xi}z|^{4+2\alpha}(1-|z|^2)^{-\alpha}}\le\\
&\int_{U^n}(1-|t|^2)^m |Df(t)|\frac{(1-|t|^2)^\alpha
\omega(1-|t|)}{(1-|t|^2)^{m-1+2+\alpha}}\, dm_{2n}(t)=\\
&\int_{U^n}\frac{|Df(t)| \omega(1-|t|)\, dm_{2n}(t)}{(1-|t|)}=\Vert
f\Vert_{B_1(\omega)}\cdot \Vert g\Vert_\infty.
\end{split}
\]
Now using again the described technique of  selection of $f_r$ by (\ref
{ef1}) for $p=1$ and $V^n$, taking $f_r(\zeta)\equiv |f_r(\zeta) |$ we
get
\[
\begin{split}
\|B_\alpha (f_r)\|_{L^1(\omega)}&\geq
\int_{U^n}\omega (1-|z|)(1-|z|^2)^{\alpha+2}
\int_{V^n}\frac{(1-|\zeta |)^\alpha}{|1-\overline\zeta z|^{2\alpha
+4}}|f_r(\zeta)|dm_{2n}(\zeta)dm_{2n}(z)\\&\geq
C_1(\alpha, \omega)\frac{(1-r)^{\alpha }}{\omega
(1-r)(1-r)^2}\int_{U^n}\frac{\omega (1-|z|)(1-|z|^2)^{\alpha +2}}{|1-r
z|^{2\alpha +4}} dm_{2n}(z).\nonumber
\end{split}
\]
As in the case of little Hankel operators, assumption of the converse
results in a contradiction.

\hfill
$\square$

\begin{corollary} Let  $\alpha_j>\alpha_{\omega_j},\,1\leq j\leq n,
\,\,g\in L^\infty (U^n).$ Then $B^\alpha_g$ is bounded on $L^1(\omega).$
\end{corollary}

                     \noindent
                     Anahit Harutyunyan \\
                     Fac. for Inf. and Appl. Math.  \\ University of Yerevan \\ Alek Manukian 1 \\
                     Yerevan 25
                     Armenia \\
                     e-mail: anahit@ysu.am
                     \vspace{4mm} \\
                      George Marinescu  \\
                     Mathematical Institute \\
                     University of Cologne \\
                     Weyertal 86 - 90\\
                     D 50931 Cologne \\
                     Germany \\
                     e-mail: gmarines@math.uni-koeln.de

\end{document}